\documentclass{article}      
\usepackage{graphicx}
\usepackage{amsmath}
\usepackage{amssymb}
\usepackage{amscd}
\usepackage{hyperref}
\usepackage{color}
\usepackage{enumerate}
\usepackage{subcaption}
\usepackage{setspace}
\usepackage{authblk}
\usepackage[top=1in, bottom=1in, left=1in, right=1in]{geometry}

\usepackage{mathptmx}

\newtheorem{remark}{Remark}

\newcommand{\eps}{\varepsilon}

\begin{document}

\title{Relaxation oscillations in an idealized ocean circulation model \thanks{This work was supported by NSF grants DMS-0940363 and DMS-1239013.}}

\author[1]{Andrew Roberts}
\author[2]{Raj Saha}

\affil[1]{Cornell University, Department of Mathematics}

\affil[2]{Bates College, Department of Physics and Astronomy}

%
		


\maketitle



\begin{abstract}

This work is motivated by a desire to understand transitions between stable equilibria observed in Stommel's 1961 thermohaline circulation model.  We adapt the model, including a forcing parameter as a dynamic slow variable.  The resulting model is a piecewise-smooth, three time-scale system.  The model is analyzed using geometric singular perturbation theory to demonstrate the existence of attracting periodic orbits.  The system is capable of producing classical relaxation oscillations as expected, but there is also a parameter regime in which the model exhibits small amplitude oscillations known as canard cycles.  Forcing the model with obliquity variations from the last 100,000 years produces oscillations that are modulated in amplitude and frequency. The output shows similarities with important features of the climate proxy data of the same period.
\end{abstract}


\section{Introduction}
Climate variability is an important aspect of the climate system.  An understanding of this variability, specifically with regard to glacial millennial climate change, has remained elusive \cite{paleoclimate}.  In different eras of Earth's history, the variations themselves change in both period and amplitude.  It is even possible to have small oscillations superimposed over larger oscillations.  Increasingly, scientists are utilizing improved technology to study the climate system through high-powered computer simulations.  Large scale oscillations and critical transitions, however, are often better understood by examining conceptual models that can be studied analytically.  Crucifix reviews key dynamical systems concepts and their applications to paleoclimate problems in \cite{crucifix}, mentioning relaxation oscillations in particular.

One of the concepts described in \cite{crucifix} is that of the relaxation oscillator, requiring a separation of time scales.  The separation of time scales gives rise to a fast/slow system that can be analyzed using {\it geometric singular perturbation theory} (GSPT).   A relaxation oscillator is characterized by ``long periods of quasi-static behavior, followed by short periods of rapid transition'' \cite{ksRO}.  Since relaxation oscillators only require a minimum of two variables (one fast, one slow), much of the full climate system might not be relevant to understanding the key mechanism(s) behind a particular oscillation.  Rather, it should be possible to consider a relevant underlying subsystem that is analytically tractable.

\begin{figure}[t]
\includegraphics[width=\textwidth]{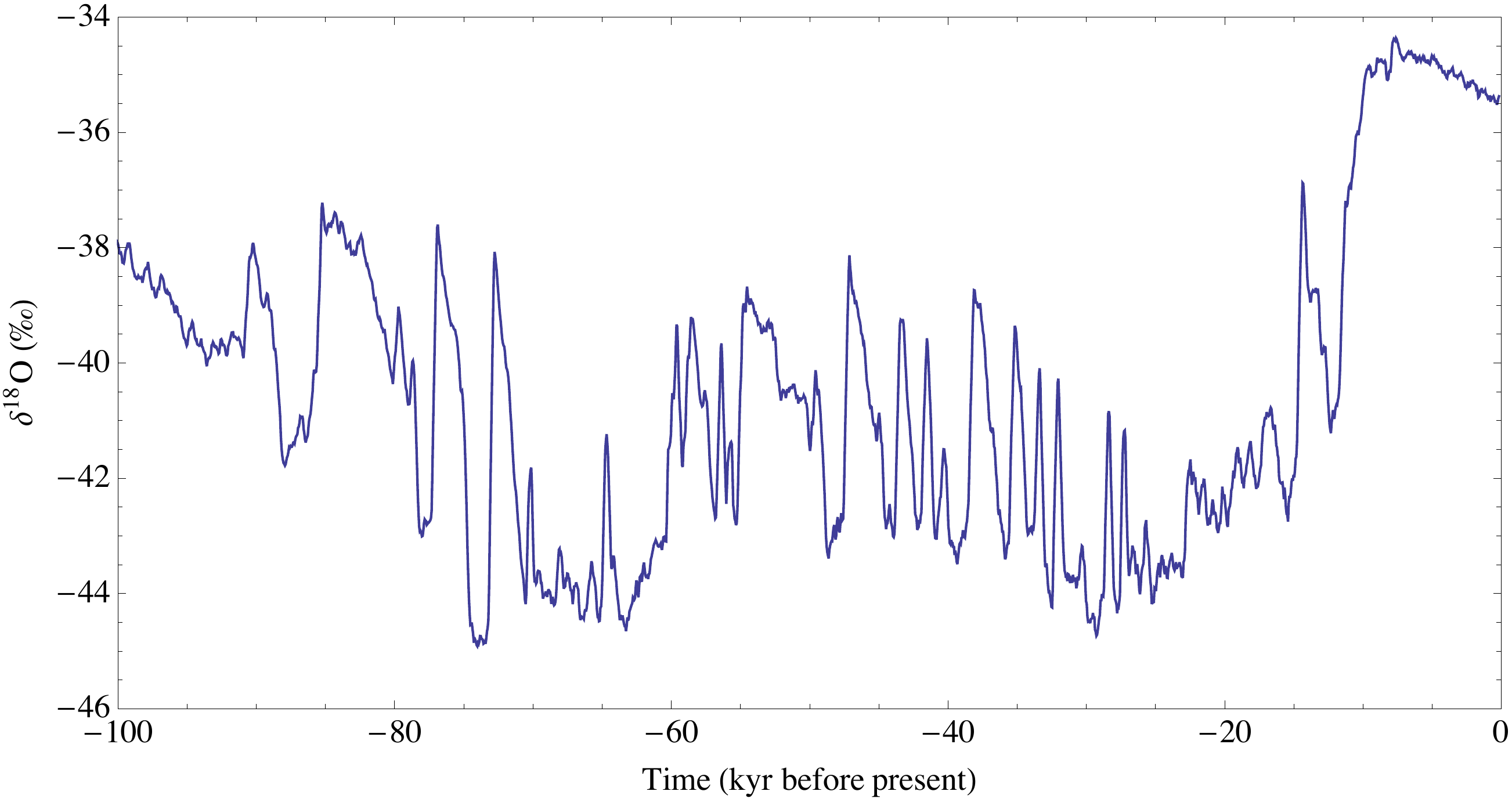}
\caption[Oxygen isotope data from Greenland (NGRIP)]{Oxygen isotope data from Greenland (NGRIP) \cite{ngrip}.}
\label{stomFig:o18}
\end{figure}

The climate record for the last 100 kiloyears (kyr) show a number of quasi-periodic pulses of abrupt warming \cite{dansgaard1993} known as Dansgaard-Oeschger (D-O), as shown in Figure \ref{stomFig:o18}. The temperature anomalies associated with D-O events were manifested with the greatest intensity in the North Atlantic region. At the onset, surface air temperatures increased sharply over a few decades followed by a gradual decline and finally a sharp drop over several centuries \cite{paleoclimate}, a characteristic that is reminiscent of relaxation oscillations. Although the timing and amplitude of individual events show considerable variations over the span of the glacial period, a fundamental pacing period of approximately 1,500 years is observed in the data \cite{schulz2002}. The lack of correlation to solar output or other astronomically forced periodicities suggest that the driver of D-O events is an internal climatic mechanism.

Since the effects of these critical climate transitions are most dramatic in the North Atlantic, scientists have hypothesized that D-O events are accompanied by changes in ocean circulation in the North Atlantic.  The bistability of the circulation in the North Atlantic was first demonstrated by Stommel in 1961 \cite{stommel}.  Physical oceanographers have provided a vast array of models capturing various mechanisms that can cause the circulation to oscillate between the two steady states in Stommel's model.  Dijkstra and Ghil surveyed many of these models in \cite{dijkstraghil}.  Some models generate oscillations as a result of intrinsic ocean dynamics \cite{dV06}.  Other models generate oscillations due to changes in freshwater forcing, be it periodic \cite{periodicfreshwater} or stochastic \cite{noisywater}.  In \cite{saltzman81}, Saltzman, Sutera, and Evenson argue that thermal effects are the driving force behind the oscillations.  Additionally, Saha shows that an ocean-ice feedback mechanism can generate oscillations  \cite{thsaha,ppsaha}.

The aims of this paper are to demonstrate the existence of attracting periodic orbits in two adaptations of Stommel's model.  In the first adaptation, we incorporate the `freshwater forcing' parameter---actually a ratio of precipitation forcing to thermal forcing---as a dynamic slow variable.  The adapted model is a three time-scale model with three variables (1 fast, 1 intermediate, and 1 slow).  We use GSP to reduce the model to a 2D fast/slow system and show that, for a certain parameter range, the reduced system has an `2'-shaped fast nullcline (a non-differentiable version of an `S'-shaped nullcline).

 Fast/slow systems with `S'-shaped nullclines---the Van der Pol system \cite{vdp} is probably the most-studied example---commonly exhibit relaxation oscillations \cite{pwlCanards,ksRO,rotstein,sw04}.  In the context of fast/slow systems, relaxation oscillations can be thought of as a series fast transitions between stable slow manifolds.  However, more complicated behavior is possible due to the existence of {\it canard} trajectories.  A canard is trajectory that passes from an attracting slow manifold to a repelling slow manifold, and remains near the repelling slow manifold for $\mathcal{O}(1)$ time \cite{mmoSurvey,ksRO,aar1}.  One reason the Van der Pol system has been a popular model for examination is the existence of periodic orbits that contain canard segments, called {\it canard cycles}.  \cite{benoitEA,pwlCanards,dumort96,eck,ksRO}.  We analyze the first adaptation of Stommel's modle to find conditions under which the model has either a relaxation oscillation or a canard cycle.

In the second adaptation of the model, we apply a 41kyr periodic forcing corresponding to changes in the Earth's obliquity.  Obliquity is not the driving force behind the D-O events, however it does play a role in modulating the amplitudes of the oscillations.  Furthermore, the forced model indicates that canard cycles might be a component and precursor of D-O events, possibly as high frequency oscillations in sea ice extent \cite{ppsaha}, although they may not be directly observable in the climate proxy record.

The relevant background material from GSPT is discussed in section 2.  Section 3 describes the physical mechanisms in Stommel's 1961 model.  We also analyze the model using geometric singular perturbation theory, following the analysis of Glendinning \cite{glendinning}.  In Section 4, we develop and analyze adapted models.  Finally, the paper concludes with further discussion in Section 5.

\section{Fast/slow dynamics}
A fast/slow system is a system of the form
\begin{equation}
	\label{gen_fast}
	\begin{array}{l}
		x' = f(x,y,\eps) \\
		y' = \eps g(x,y,\eps),
	\end{array}
\end{equation}
where the prime $'$ denotes $d/dt$.  Here, $x \in \mathbb{R}^n$ is a vector of fast variables and $y \in \mathbb{R}^m$ is a vector of slow variables.  Typically, $f$ and $g$ are assumed to be smooth functions.  If time is rescaled by $\eps \neq 0$, i.e. $\tau = \eps t$, then \eqref{gen_fast} becomes 
\begin{equation}
	\label{gen_slow}
	\begin{array}{l}
		\eps \dot{x} = f(x,y,\eps) \\
		\dot{y} = g(x,y,\eps),
	\end{array}
\end{equation}
where the dot $\dot{ \ \ }$ denotes $d/d\tau$.  The fast system \eqref{gen_fast} and the slow system \eqref{gen_slow} are equivalent as long as $\eps \neq 0$.  However, insight can be gained by looking at the $\eps = 0$ limit.  Indeed, the reduction in phase space dimensions from considering the limit may turn an analytically intractable problem into a tractable one.

As $\eps \rightarrow 0$, \eqref{gen_fast} approaches 
\begin{equation}
	\label{fast_0}
	\begin{array}{l}
		x' = f(x,y,0) \\
		y' = 0.
	\end{array}
\end{equation}
which is called the {\it layer problem.}  The slow system \eqref{gen_slow} approaches
\begin{equation}
	\label{slow_0}
	\begin{array}{l}
		0 = f(x,y,0) \\
		\dot{y} = g(x,y,0),
	\end{array}
\end{equation}
which is called the {\it reduced problem.}  In both cases, the set $M_0 = \{ f(x,y,0) = 0 \}$---called the {\it critical manifold}---is special.  $M_0$ is the set of critical points of \eqref{fast_0} and the algebraic set on which \eqref{slow_0} is defined.  One can obtain a caricature of \eqref{gen_fast} or \eqref{gen_slow} by allowing the layer problem to equilibrate and then considering the reduced problem.  The theory of GSP describes the manner in which the two subsystems are pieced together provided 
\begin{equation}
	\label{normHyp}
	\det \left( \frac{ \partial f}{\partial x} |_{M_0} \right) \neq 0.
\end{equation}
 A survey of the basic theory of GSP can be found in \cite{gsp}.
 
\subsection{Relaxation Oscillations and Canards}
Since $M_0$ is a set of critical points of the layer problem, we can discuss its stability in a natural way.  The condition \eqref{normHyp} is a non-degeneracy condition that says $M_0$ is hyperbolic with respect to the fast dynamics.  $M_0$ is attracting (resp. repelling) where its points are attracting (resp. repelling) critical points of the layer problem.  The interesting and scientifically relevant mathematics often occur due to degenerate points, and that is the case in the models we will examine in this paper.  

Perhaps the most common type of degeneracy is a folded critical manifold.  A {\it fold} is the set of saddle-node bifurcations of the layer problems.  At the fold $f_x$ is singular, and it coincides with a co-dimension one set where $M_0$ is attracting on one side and repelling on the other.  If $M_0$ has two folds, we say it is `S'-shaped.  In planar systems, `S'-shaped critical manifolds can produce relaxation oscillations; in higher dimensional systems, even more complicated oscillatory behavior is possible.  The relaxation oscillations result from a Hopf bifurcation that occurs as the slow nullcline passes through a fold point.  If the Hopf bifurcation is supercritical, the initial orbits born of the bifurcation are small amplitude orbits called canard cycles.  Over an exponentially small parameter range, the canard cycles will grow into large relaxation oscillations.  This rapid transition from small orbits to relaxation oscillations is a phenomenon know as a {\it canard explosion} \cite{ksRO}.  If the Hopf bifurcation is subcritical, the attracting periodic orbits created will be relaxation oscillations.  Additionally, in higher dimensions generic fold points can produce relaxation oscillations without a canard explosion.

The models we analyze in this paper do not have smooth vector fields due to an absolute value term.  The set where the vector field is not differentiable is called the {\it splitting line}.  As we will see in our models, it is possible for a `corner point' of the critical manifold to occur when it intersects the splitting line.  In \cite{aar1}, Roberts and Glendinning describe the possible oscillatory behavior in a fast/slow system with a piecewise-defined `2'-shaped critical manifold.  That is, they consider a system of the form
\begin{equation}
	\label{generalStom}
	\begin{array}{l}
		\dot{x} = -y + F(x) \\
		\dot{y} = \eps ( x - \lambda)
	\end{array}
\end{equation}
where $$ F(x) = \left\{ \begin{array}{ll}
	g(x) & x\leq 0 \\
	h(x) & x \geq 0
	\end{array} \right. $$ 
with $g,h \in C^k, \ k \geq 1$, $g(0) = h(0) =0$, $g'(0) < 0$ and $h'(0) > 0$, and assume that $h$ has a maximum at $x_M > 0$.  The critical manifold 
$$ M_0 = \{ y = F(x) \} $$ is `2'-shaped with a smooth fold at $x_M$ and a corner along the splitting line $x = 0$. 

\begin{figure}[t!]
		\centering
		\begin{subfigure}[t]{0.3\textwidth}
			\includegraphics[width=\textwidth]{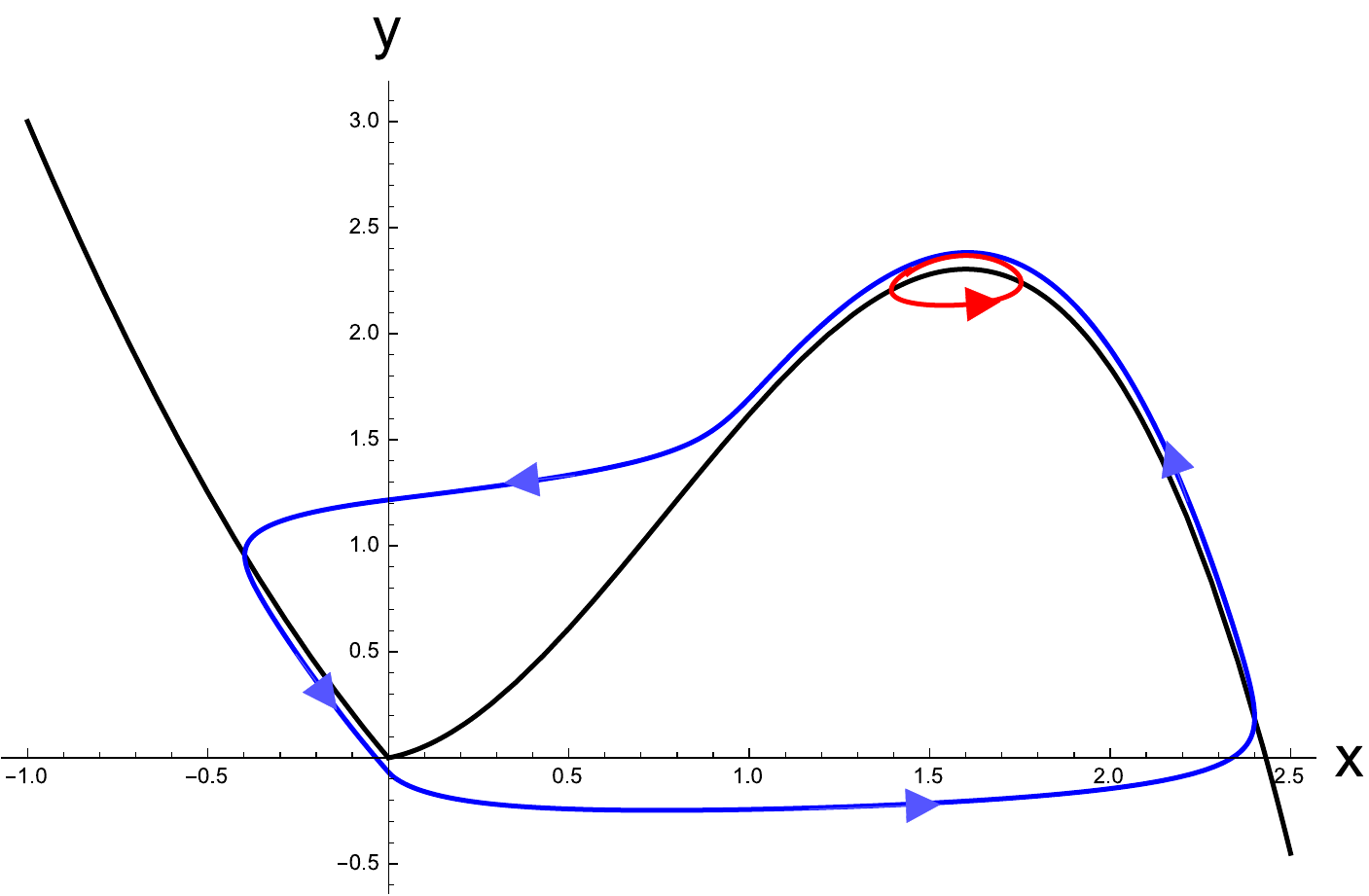}
			\caption{Examples of periodic canard trajectories resulting from a smooth Hopf bifurcation when $\eps=0.1$.  For the small cycle (red), $\lambda$ is 0.3658835 away from the critical value; for the large cycle (blue), $\lambda$ is 0.3658840 away from the critical value.}
			\label{stomFig:fold}
		\end{subfigure}
		~
		\begin{subfigure}[t]{0.3\textwidth}
			\includegraphics[width=\textwidth]{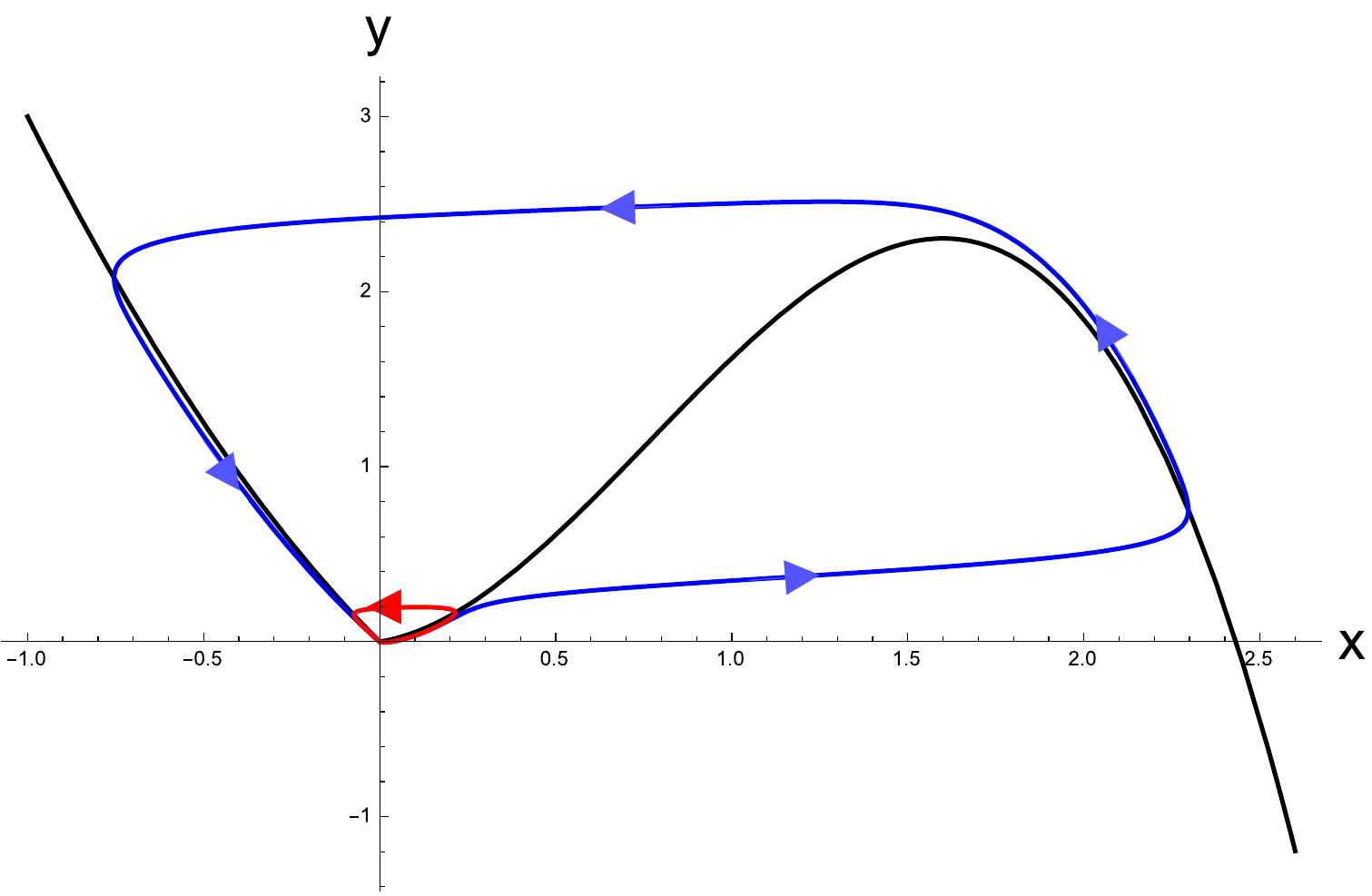}
			\caption{Examples of periodic canard trajectories resulting a nonsmooth Hopf bifurcation when $\eps=0.1$.  For the small cycle (red), $\lambda = 0.01293$; for the large cycle (blue) $\lambda=0.01295$.}
			\label{stomFig:corner}
		\end{subfigure}
				~
		\begin{subfigure}[t]{0.3\textwidth}
			\includegraphics[width=\textwidth]{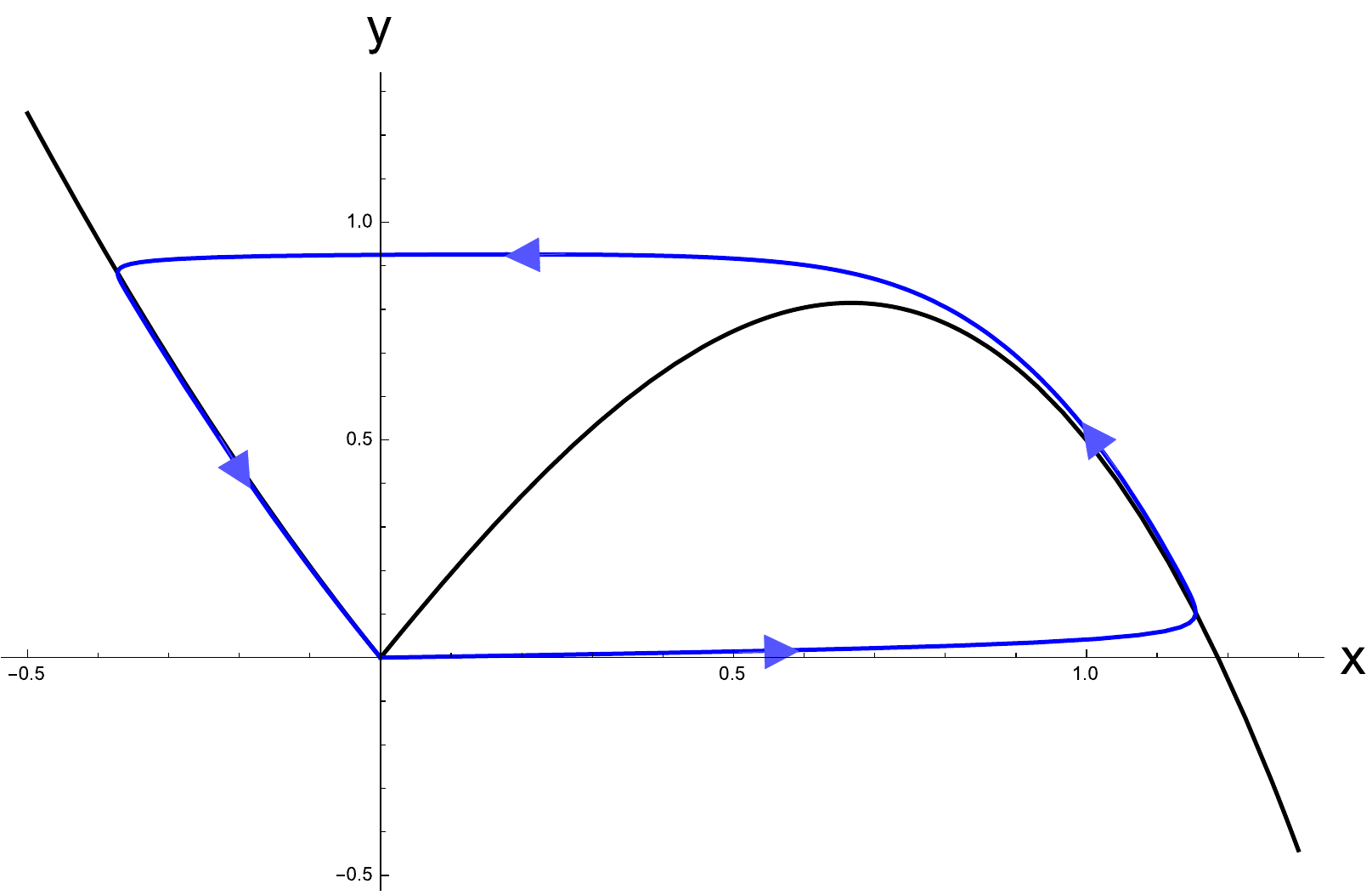}
			\caption{The stable orbit of a super-explosion (blue) with $y$-nullcline (red) when $g(x) = (x-1)^2-1,$ $h(x) = -(x+1)^2 (x-1.5)-1.5,$ $\eps=0.01$, and $\lambda=10^{-7}.$}
			\label{stomFig:superExp}
		\end{subfigure}
		\caption{Example of canard-related oscillations resulting from supercritical bifurcations in \eqref{generalStom}.  Figures created using NDSolve in Mathematica.}
		\label{stomFig:nsCanards}
	\end{figure}

The bifurcation that occurs when $\lambda = x_M$ is a traditional Hopf-bifurcation, and any canard explosion phenomenon is essentially classical \cite{aar1}.  The bifurcation that occurs when $\lambda = 0$ is not a classical bifurcation since the vector field is not $C^1$.  At this nonsmooth bifurcation, the local geometry near the origin is extremely important in determining the nature of the periodic formed as $\lambda$ increases through 0.  When $\lambda <0$ (but small), the system has a globally attracting fixed point, and depending on $g'(0)$, the point will either be a stable node (2 real negative eigenvalues) or a stable focus (2 stable complex eigenvalues).  Similarly, when $\lambda > 0$ (but small), the equilibrium point will either be an unstable node (2 real positive eigenvalues) or unstable focus (2 unstable complex eigenvalues).  However, since the vector field is not differentiable, the real part of the eigenvalues never passes through 0, allowing for bifurcations that are different from those found in smooth systems.  Figure \ref{stomFig:nsCanards} depicts the possible canard-related orbits in systems like \eqref{generalStom}, described below.

In essence, the bifurcation can be determined by the type of equilibrium that exists when $\lambda > 0$, and the criticality is determined by the equilibrium when $\lambda <0$.  If an unstable focus is created, the bifurcation behaves like a nonsmooth version of a Hopf bifurcation due to the rotation that happens near the unstable equilibrium.  In case where there is a transition from a stable focus to an unstable focus, the criticality is determined by the values of $|g'(0)|$ and $|h'(0)|$, with the bifurcation being subcritical when $|g'(0)| < |h'(0)|$ and supercritical if the inequality is reversed.  If the bifurcation is described by a transition from a stable node to an unstable focus, the bifurcation will always be supercritical since there is insufficient rotation near the equilibrium to produce a periodic orbit before bifurcation.  If the nonsmooth Hopf bifurcation is supercritical, the bifurcation will create small amplitude periodic orbits that rapidly grow in amplitude as $\lambda$ increases away from 0, similar to the classical canard explosion phenomenon.

If the bifurcation creates an unstable node, then the system will transition immediately into a relaxation oscillation regime in a bifurcation that is called a {\it super-explosion}.  The term super-explosion refers to the fact that the bifurcation forgoes the classical canard explosion.  If the bifurcation transitions from a stable node to an unstable node, then the bifurcation will be supercritical (i.e., the relaxation periodic orbit only exists after the bifurcation).  However, the the bifurcation transitions from a stable focus to an unstable node, the rotation before bifurcation causes the bifurcation to be subcritical, and a stable relaxation periodic orbit coexists with the stable equilibrium for some $\lambda < 0$ but small.

It is worth noting that the nature of the equilibrium points, and hence the nature of the bifurcation, is highly dependent upon $\eps$.  In contrast to classical canard phenomena, $\eps$ can be too small to create canard cycles in the nonsmooth case.  In particular, canard cycles require the inequality $|h'(0)| < 2 \sqrt{\eps}$ to hold.  For a more in-depth discussion of these bifurcations we direct the reader to \cite{aar4,aar1}.

\begin{remark}
This nonsmooth bifurcation has been generalized to allow for more complicated slow dynamics than occur in \eqref{generalStom}, possibly moving the bifurcation location of the bifurcation point.  The more general theorems can be found in \cite{aar4}, however the versions presented above are sufficient for the models we consider here.
\end{remark}

The remainder of the paper will be devoted to developing and analyzing a fast/slow system modeling large-scale ocean circulation.

\section{Stommel's Model}
\begin{figure}[t]
\includegraphics[width=\textwidth]{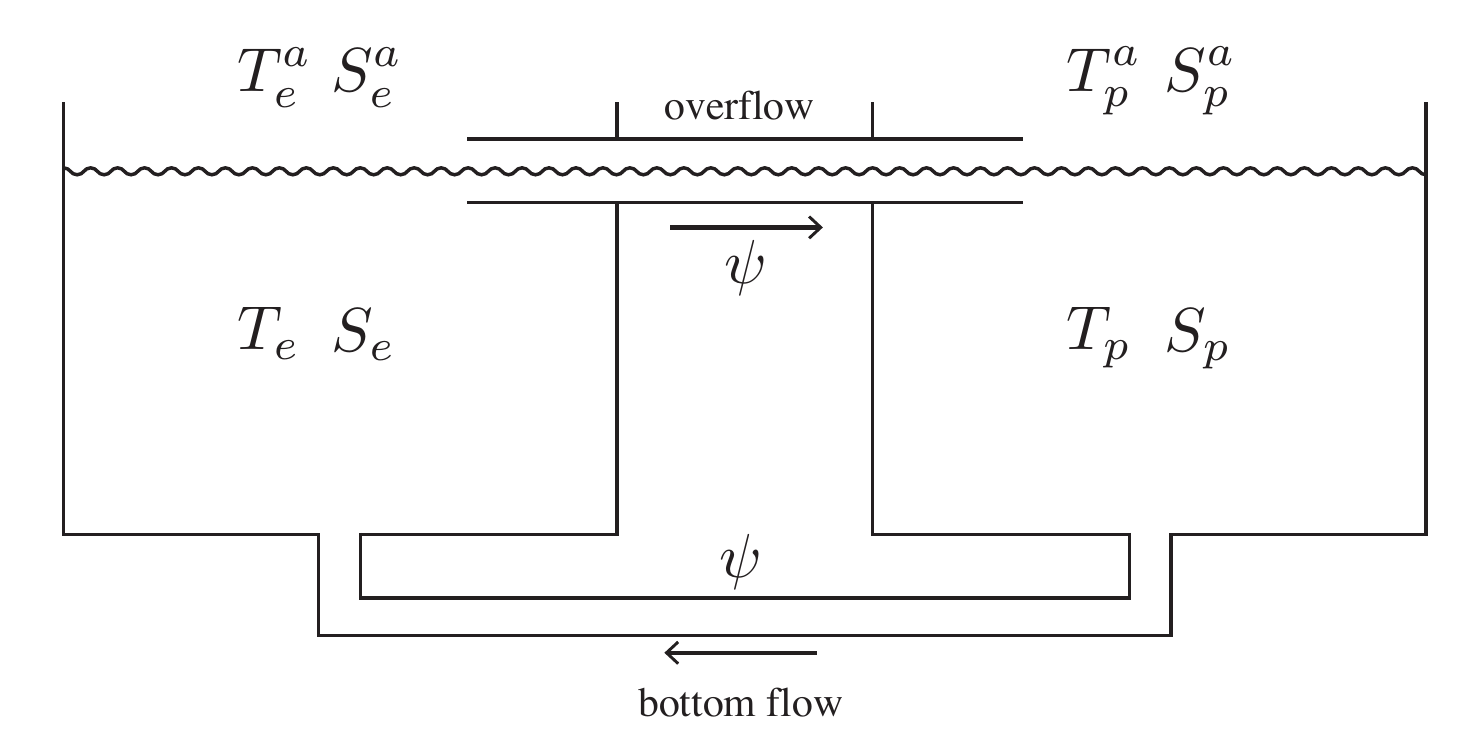}
\caption{Schematic of Stommel's model.}
\label{stomFig:schematic}
\end{figure}

	Investigating changes in ocean circulation begins with Stommel's 2-box model.  Stommel modeled the North Atlantic by partitioning it into an equatorial and a polar region.  He assumed that water near the equator would become warmer and saltier due to its interaction with the atmosphere.  Water near the pole would lose its heat to the atmosphere and have its salt concentration diluted by incoming freshwater.  Oceanic circulation causes the water in the two regions to mix, preventing either box from equilibrating with its surrounding environment.  The density difference between the boxes drives the circulation, and it is the salinity and temperature of each box that determines the density.  A schematic of the model is shown in Figure \ref{stomFig:schematic}.  The analysis in this section follows \cite{glendinning} and \cite{kg11}.  It is included here as a reminder of the mechanics of Stommel's model and a means of setting up the equations for the main results of this paper.

The equations describing the model are
	\begin{equation}
		\label{stom4}
		\begin{array}{c}
			\frac{d}{dt}T_e = R_T (T_e^a - T_e) +| \psi | (T_p - T_e) \\
			\frac{d}{dt}T_p = R_T (T_p^a - T_p) +| \psi |  (T_e - T_p) \\
			\frac{d}{dt}S_e = R_S (S_e^a - S_e) +| \psi |  (S_p - S_e) \\
			\frac{d}{dt}S_p = R_S (S_p^a - S_p) +| \psi | (S_e - S_p).
		\end{array}
	\end{equation}
Here, $T$'s are temperatures and $S$'s are salinities.  The subscripts $e$ and $p$ denote the box at the equator and pole, respectively, while the superscript $a$ denotes an atmospheric forcing term.  The strength of the circulation is given by $| \psi |$, where
	$$  	\psi = \psi_0 \left( \frac{ \rho_p - \rho_e } {\rho_0} \right).  $$
	
	The density of box $i$ is denoted $\rho_i$, and it is calculated using a linear equation of state given by
	$$  \rho_i = \rho_0[1- \alpha(T_i - T_0) + \beta(S_i - S_0)]. $$
where $\alpha$ and $\beta$ are the coefficients of thermal expansion and haline contraction respectively, and $T_0$, $S_0$ and $\rho_0$ are reference values.  	
	
The advective term $\psi$ can thus be computed to be
	\begin{equation}
		\label{circ}
			\psi = \psi_0 [ \alpha (T_e - T_p) - \beta (S_e - S_p) ].
	\end{equation}

Next, the system is reduced to have only two degrees of freedom by looking at the temperature and salinity differences between the boxes.  Defining
	\[ 
		\begin{array}{ccc}
			T = T_e - T_p 	,		& 		\hspace{1in} 	&	S = S_e - S_p, \\
			T^a = T_e^a - T_p^a, 	&		\hspace{1in} 	&	S^a = S_e^a - S_p^a,\\
			X = T_e + T_p			&		\hspace{1in}	&	Y = S_e + S_p,
		\end{array}
	\]
we see that the the equations for $dX/dt$ and $dY/dt$ decouple from the others.  $X$ and $Y$ can be found trivially, by direct integration, from their governing equations.  Therefore, the interesting dynamics of \eqref{stom4} are captured by the planar system:
	\begin{equation}
		\label{stom2}
		\begin{array}{c}
			\frac{d}{dt}T = R_T (T^a - T) - 2 | \psi | T \\
			\frac{d}{dt} S = R_S (S^a - S) - 2 |\psi | S,
		\end{array}
	\end{equation}
and \eqref{circ} becomes $$ \psi = \psi_0 (\alpha T - \beta S). $$

\begin{figure}[t]
		\centering
		\begin{subfigure}[t]{0.45\textwidth}
			\includegraphics[width=\textwidth]{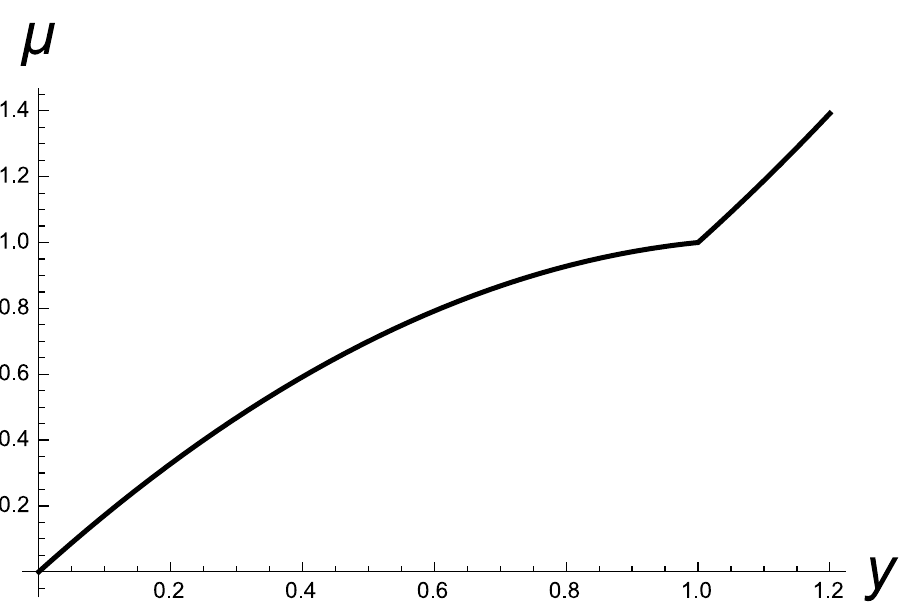}
			\caption{$A<1$}
			\label{stomFig:smallA}
		\end{subfigure}
		~
		\begin{subfigure}[t]{0.45\textwidth}
			\includegraphics[width=\textwidth]{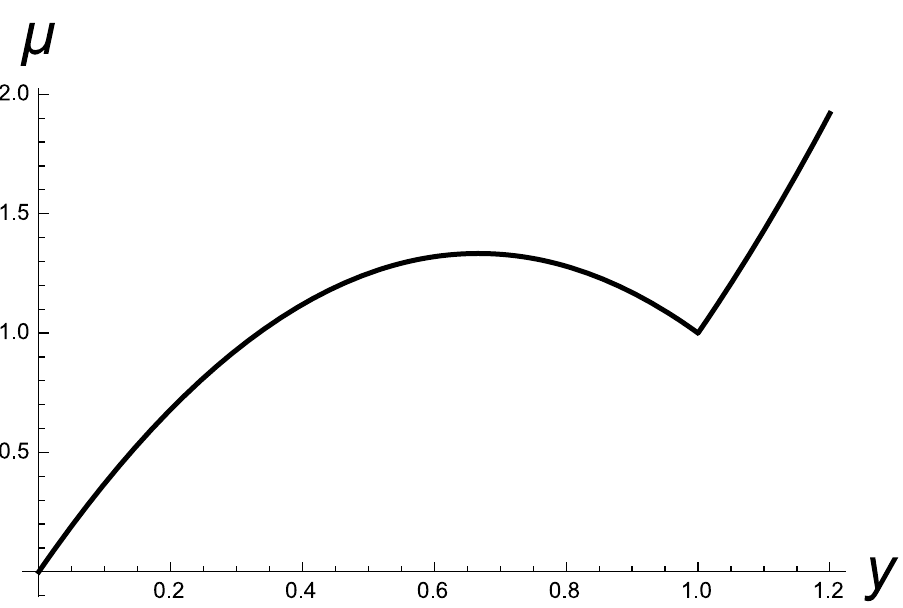}
			\caption{$A>1$}
			\label{stomFig:bigA}
		\end{subfigure}
		\caption{Graphs of \eqref{crit1} for (a) $A<1$ and (b) $A>1$.}
		\label{stomFig:As}
	\end{figure}

To non-dimensionalize the system, set
	\[
		\begin{array}{ccccc}
			x = \displaystyle{ \frac{T}{T^a}, } & \displaystyle{ y = \frac{\beta S}{\alpha T^a}, } &  \displaystyle{ \tau = R_S t ,}
			& \displaystyle{ \mu = \frac{\beta S^a}{\alpha T^a}, } &  \displaystyle{  A = \frac{ 2 \psi_0 \alpha T^a}{R_S}. }
		\end{array}  
	\]
Then the system \eqref{stom2}  becomes
	\begin{equation}
		\label{non}
		\begin{array}{l}
			\eps \dot{x} = 1 - x - \eps A | x - y | x \\
			\dot{y} = \mu - y - A | x - y | y,
		\end{array}
	\end{equation}
where $$ \eps = \frac{R_S}{R_T} \ll 1$$ is a small parameter and the dot ( $\dot{ \ \  }$ ) denotes differentiation with respect to $\tau$.  The system \eqref{non} is a fast/slow system set up to be analyzed using GSP.  In the limit as $\eps \rightarrow 0,$ $\{ x = 1 \}$ is a globally attracting, and therefore normally hyperbolic, critical manifold.  The reduced problem has one degree of freedom, so the dynamics are entirely characterized by equilibria.  The flow on the critical manifold $\{ x= 1 \}$ is given by
	\begin{equation}
		\label{red1}
		\dot{y} = \mu - y - A | 1- y | y.
	\end{equation}
Critical points occur at 
	\begin{equation}
		\label{crit1}
		\mu = \left\{ \begin{array}{ll}
			(1 + A) y - A y^2 & \text{for } y < 1 \\
			(1 - A) y + A y^2 & \text{for } y > 1
			\end{array} \right.
	\end{equation}
and the nature of the system depends on $A$, as seen in Figure \ref{stomFig:As}. 
Taking a derivative gives
	\[ 
		\frac{d\mu}{dy} = \left\{ \begin{array}{ll}
			(1 + A) - 2 A y & \text{for } y < 1 \\
			(1 - A) + 2 A y & \text{for } y > 1.
			\end{array} \right.
	\]
	If $A<1$, then the curve of equilibria $\mu = \mu(y)$ is monotone increasing.  The system \eqref{red1}, and consequently \eqref{non}, has a unique equilibrium solution.  The equilibrium is globally attracting, and it is important to remember that the solution corresponds to a unique stable circulation state (i.e., direction and strength).

However, if $A >1$ the system exhibits bistability for a range of $\mu$ values.  While $\mu(y)$ is still monotone increasing for $y > 1$, the curve has a local maximum at $y =(1+A)/(2A) < 1$.  Thus for $1 < \mu < (1+A)^2/(4A) $ there are three equilibria.  The system is bistable with the outer two equilibria being stable, and the middle equilibrium being unstable.  The sign of $\psi$---determined by the relative strengths of the thermal and salinity gradients, $x$ and $y$ respectively in the dimensionless system---differs at each of the two stable states.  The stable equilibrium for $\psi < 0$ is called the haline state, since the circulation is driven by a stronger salinity gradient.  When $\psi > 0$, the circulation is driven by a stronger temperature gradient, and the system is in the thermal state.  In the bistable regime, there is a stable thermal state as well as an unstable thermal state.  As mentioned in the introduction, some oceanographers attempt to explain oscillations using only oceanic processes, however Figure \ref{bifurcation} suggests $\mu$ is the key to generating such oscillations.  

\section{Results}
	\begin{figure}[t]
		\begin{center}
			\includegraphics[width=0.5\textwidth]{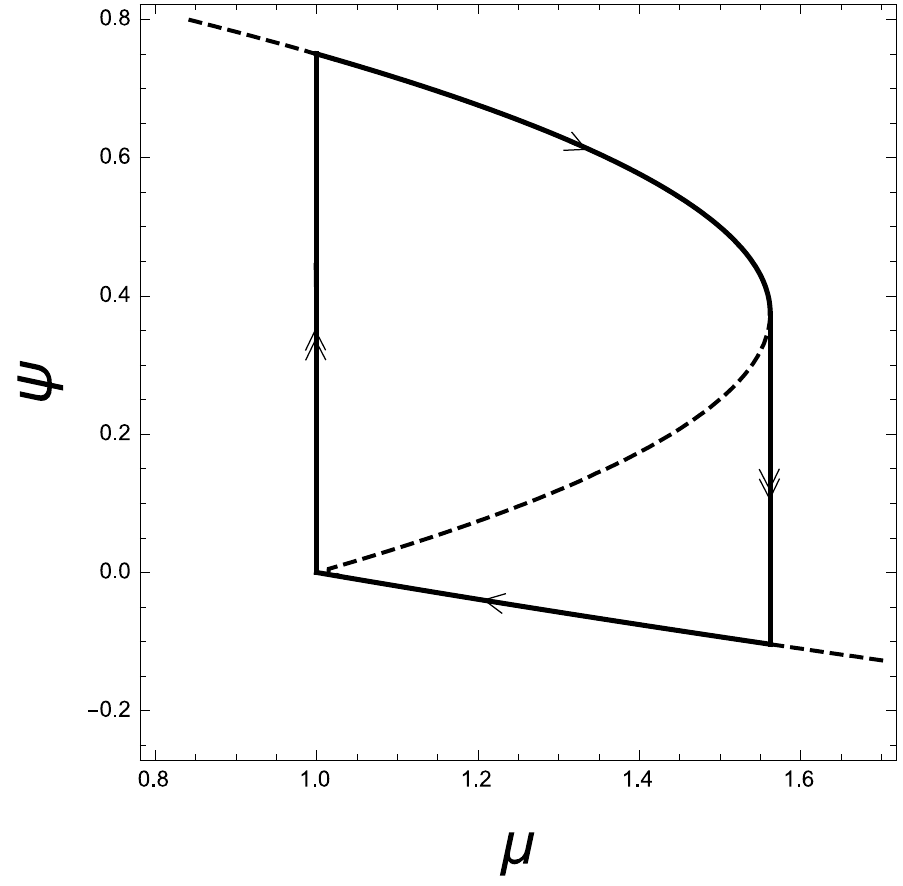}
		\caption{Bifurcation diagram for \eqref{non} (dashed), with a hysteresis loop (solid black) overlay.  $\psi = \alpha \psi_0(x-y).$}
		\label{bifurcation}
		\end{center}
	\end{figure}
If $\mu$ is the key to oscillations, there may be an intrinsic feedback mechanism that causes $\mu$ to change.  Recall that $\mu$ is the ratio of the effect of atmospheric salinity forcing on density to that of atmospheric temperature forcing on density.  The idea is to consider how the state of the ocean affects its interaction with the atmosphere.  Typically, in a coupled ocean-atmosphere model, the atmosphere is the fast component and the ocean is the slow component (see \cite{vanveen}, for example).  However, in Stommel's model, $\mu$ is considered constant.  Therefore, if a model is going to incorporate $\mu$ as a dynamic variable, it should vary on a slower time scale than the other variables in the model.  The physical intuition is to consider the variation of long term average behavior in the atmosphere.   This is reminiscent of the idea of a {\it dynamic bifurcation}, where the bifurcation parameter varies slowly according to a prescribed equation \cite{ma2005}.  Here, we assume that variations in $\mu$ will depend on the state variables $x,y$.  A general system of this form is 
	\begin{equation}
		\label{general3}
		\begin{array}{l}
			 x' = 1 - x - \eps A | x - y | x \\
			 y' = \eps( \mu - y - A | x - y | y )\\
			\mu' =  \eps \delta f(x,y,\mu,\delta,\eps),
		\end{array}
	\end{equation}
where $\delta \ll 1$ is another small parameter.  This system is a three time-scale model where $x$ is fast, $y$ is intermediate, and $\mu$ is slow.  If $f(x,y,\mu,0,0)$ behaves in a desirable manner, one can turn the hysteresis loop implied by Figure \ref{bifurcation} into a dynamic periodic orbit. The important question to answer is ``how do variations in $x$ and $y$ affect $\mu$?''. One way in which the state of the atmosphere can be affected by the state of the ocean is with the presence of sea ice. A cooling of the ocean produces sea ice which on a short time scale insulates the ocean from exchanging heat with the atmosphere. Over longer time scales the net effect of multiple cycles of sea ice growth and melt could result in an overall heat loss from the ocean and enhanced convective mixing \cite{ppsaha}. The ventilated heat from the ocean would warm the polar atmosphere and decrease the atmospheric temperature gradient, $T^a$. Recalling that 
	$$ \mu = \frac{\beta S^a}{\alpha T^a}, $$
it can be seen that an increase in $x$, the scaled temperature gradient of the ocean due to cooling of the polar ocean would increase $\mu$.

A strong salinity gradient in the ocean acts as a negative feedback to the atmospheric salinity (freshwater) gradient as saltier water tends to be less evaporative. Thus an increase in $y$, the scaled oceanic salinity gradient, would decrease $\mu$. Both of the feedback mechanisms described above could thus produce the dependence of $f$ on $x$ and $y$ as required for oscillations, i.e., 
	\begin{equation}
		\label{effects}
		\begin{array}{ccc}
			\displaystyle{ \frac{\partial f}{\partial x} > 0} & \text{and} & \displaystyle{ \frac{\partial f}{\partial y} < 0}.
		\end{array}
	\end{equation}
	
If this condition is implemented in the simplest possible way, then there is a parameter regime in which the system has a unique periodic orbit.  Taking $f$ to be a linear function of $x$ and $y$, \eqref{general3} becomes
	\begin{equation}
		\label{lin3}
		\begin{array}{c}
			x' = 1 - x - \eps A | x - y | x \\
			y' = \eps( \mu - y - A | x - y | y) \\
			\mu' = \eps \delta (1+ ax - by).
		\end{array}
	\end{equation} 	
	
As in the previous section, $\{ x = 1\}$ is still an attracting critical manifold.  However, the reduced problem,
	\begin{equation}
		\label{red2}
		\begin{array}{c}
			\dot{y} =  \mu - y - A | 1 - y | y \\
			\dot{\mu} =  \delta (1+ a - by),
		\end{array}
	\end{equation}
is now itself a fast/slow system which is analyzed using GSP.

The critical manifold of \eqref{red2} is given by
	\begin{equation}
		\label{m0}
		M_0 = \{ \mu = y + A | 1 - y | y \}.
	\end{equation}
We point out that this is precisely the curve of equilibria from \eqref{non} in the previous section (depicted in Figure \ref{stomFig:As}), and this relationship allows for easy comparison with the previous section.  First, the shape of $M_0$ depends on the parameter $A$ in the same way as the curves of equilibria in Figure \ref{stomFig:As}.  Second, stable (respectively unstable) branches of $M_0$ correspond to stable (respectively unstable) segments of the curve of equilibria from \eqref{non}.  Thus, if $A<1$, $M_0$ is everywhere attracting; if $A>1$, $M_0$ is bistable, with the outside branches being attracting, and the middle branch (the segment between the fold and the corner in Figure \ref{stomFig:bigA}) is repelling.

The dynamics on $M_0$ depend only on parameters and $y$, which is fast in the reduced problem.  Therefore, the key to resolving the slow flow is the location of the $\mu$ nullcline.  The cases $A<1$ and $A>1$ are treated separately, due to the different shape of the critical manifold as depicted in Figure \ref{stomFig:As}.

\subsection{ Globally Attracting Critical Manifold---Equilibration Regime}

Here we consider the case where $A<1$.  Since the only two branches of $M_0$ are both attracting for $A<1$, one expects that system \eqref{red2} should behave as a 1-dimensional system.  The intersection of the $\mu$ nullcline and the critical manifold should be a globally attracting critical point, however GSP cannot be applied ``out of the box'' due to the non-differentiability of the vector field.  Instead, we consider two distinct smooth dynamical systems: (1) where $|1-y| = 1-y$  and (2) where $|1-y| = y-1$.  The system where $|1-y| = 1-y$ agrees with \eqref{red2} when $y<1$.  Similarly the system where $|1-y| = y-1$ agrees with \eqref{red2} when $y>1$.  

The dynamics of  \eqref{red2} is obtained by taking trajectories from the two smooth systems and cutting them along the line $y=1$---called the splitting line---where both smooth systems agree with the system of interest.  We then paste the relevant pieces together along the splitting line.  Since the right-hand side in \eqref{red2} is Lipschitz, trajectories pass through the splitting line in a well-defined manner due to uniqueness of solutions.  
	\begin{figure}[t!]
		\centering
		\begin{subfigure}[t]{0.45\textwidth}
			\includegraphics[width=\textwidth]{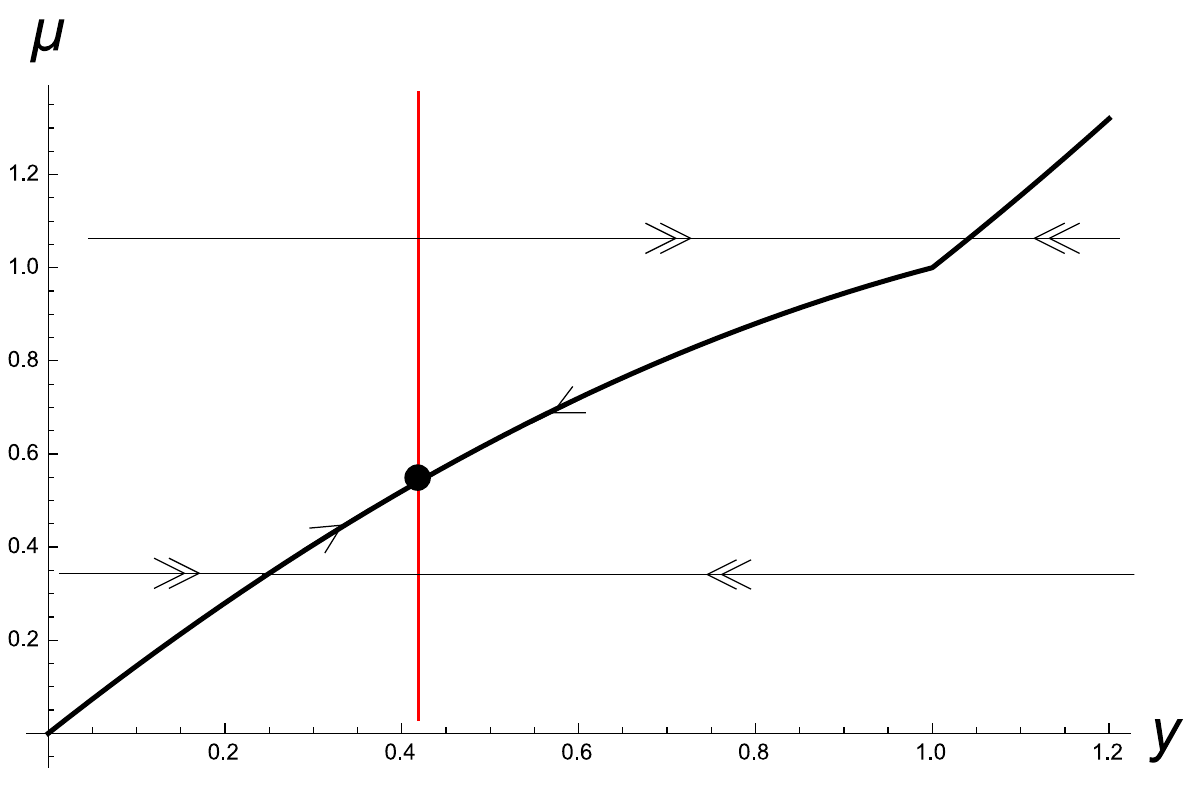}
			\caption{$\displaystyle{\frac{1+a}{b}<1}$}
			\label{stomFig:smallA_low}
		\end{subfigure}
		~
		\begin{subfigure}[t]{0.45\textwidth}
			\includegraphics[width=\textwidth]{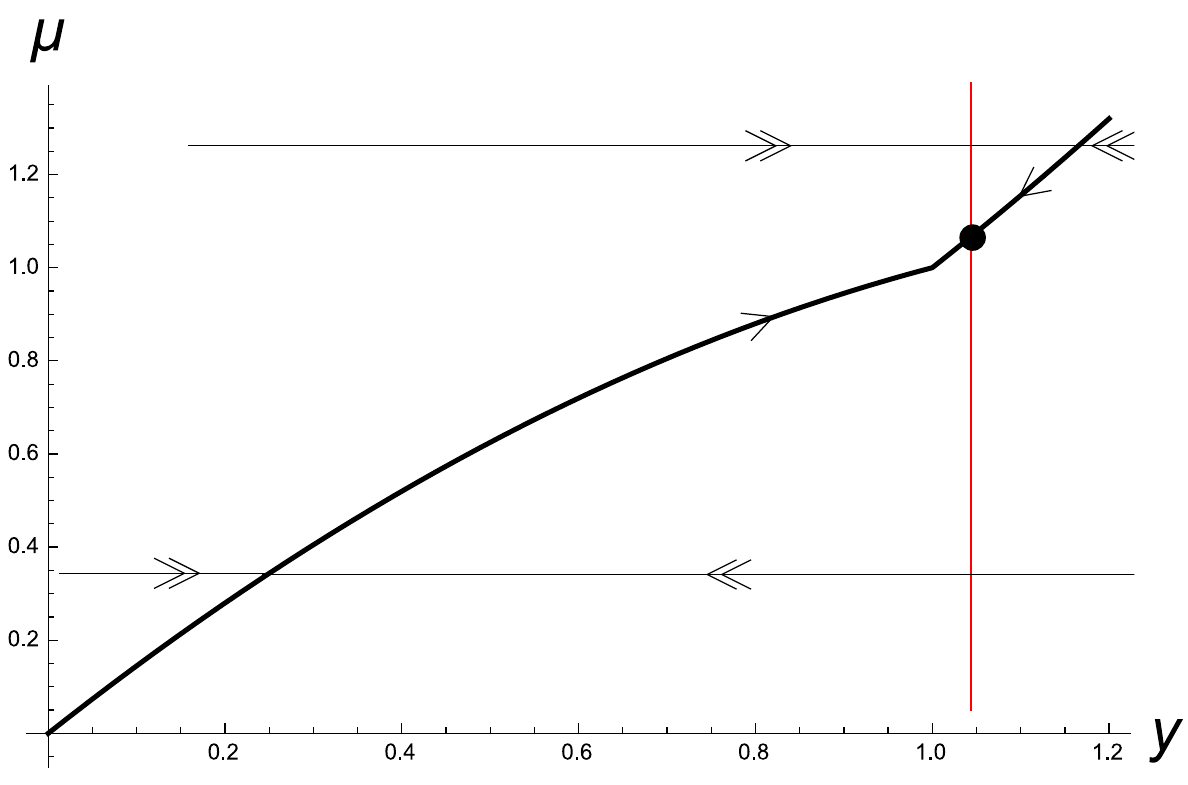}
			\caption{$\displaystyle{\frac{1+a}{b}>1}$}
			\label{stomFig:smallA_high}
		\end{subfigure}
		\caption[Possible phase spaces of \eqref{red2} for $A<1$ and $1+a \neq b$]{Possible phase spaces of \eqref{red2} for $A<1$ and $1+a \neq b$.  The red line is the $\mu$ nullcline.  The black arrows indicate fast dynamics, and the blue arrows indicate slow dynamics.}
		\label{stomFig:As2}
	\end{figure}
	
Now, GSP can be applied to both of the smooth systems, which will produce two critical manifolds that intersect (in the singular limit) when $y=1$.  $M_0$ defined in \eqref{m0} is obtained by taking the relevant critical manifold on either side of the splitting line.  The $\mu$ nullcline is the vertical line
$$y = \frac{1+a}{b}.$$
We see that $\mu$ is increasing to the left of this line, and decreasing to the right of this line as in Figure \ref{stomFig:As2}.  Therefore the system has a globally attracting equilibrium if $A<1$.

As in the previous section, the fast/slow analysis essentially determines that the ocean temperature gradient rests near the atmospheric temperature gradient.  A consequence of the equilibration regime is that the ocean salinity gradient also reaches equilibrium.  Recall that the circulation variable, $\psi = \alpha \psi_0 (x - y) \approx \alpha \psi_0 (1 - y),$ is determined by the ocean temperature and salinity gradients.  An important physical implication of the equilibration regime is that the circulation has fixed strength of approximately $\alpha \psi_0 [1-(1+a)b^{-1}].$

\subsection{Bistable Critical Manifold---Oscillating Regime}
If $A > 1$, the system is much more interesting due to the `2'-shaped critical manifold.  To simplify the analysis, we rewrite \eqref{red2} as 
\begin{equation}
	\label{canards}
	\begin{array}{rl}
		\dot{y} =& \mu - y - A |1-y| y \\
		\dot{\mu} =& \delta_0 ( \lambda - y),
	\end{array}
\end{equation}
where $\delta_0 = \delta b$ and $\lambda = (1+a)/b.$  

The simple coordinate change $(\bar{y}, \bar{\mu}) = (y-1, - \mu)$ transforms \eqref{canards} into the form of \eqref{generalStom}, however, for the purpose of physical interpretation, we will analyze \eqref{canards}.  Ultimately, we seek to understand the long time behavior of \eqref{canards} under the assumptions that $A>1$, $0 < \delta_0 \ll 1$, and $\lambda > 0$ is fixed.  We will see that, depending on the value of $\lambda$ the system can be in either an equilibration regime or an oscillating regime.

Define $F_{\pm}(y) = y \pm A (1-y) y$, and 
$$ F(y) = \left\{ \begin{array}{ll}
	F_+(y) & y< 1 \\
	F_- (y) & y  > 1 \end{array} \right. $$
Then  \eqref{canards} always has a unique equilibrium at $(y_0,\mu_0) = (\lambda, F(\lambda) )$.  Direct computation shows that the Jacobian of \eqref{canards} is 
  \begin{equation} 
  	\label{jac}
  	J ( \lambda, F(\lambda) ) = \left( \begin{array}{cc}
	 -F'(\lambda) & 1 \\
	  -\delta_0 & 0
	\end{array} \right),
  \end{equation}
We see that $\det{J} > 0$ everywhere, and $\text{Tr}(J) < 0$ when $F'(\lambda) > 0.$  Therefore, for $\lambda > 1,$ we have an attracting equilibrium.  When  $\lambda =1$ the equilibrium is attracting from the right (haline state), but repelling from the left (thermal state).  All trajectories in the thermal state will eventually be returned to the haline state above the $y$ nullcline.  Since $F'_-(1)  > 1$, we can assume $[ F'_- (1)]^2 > 4 \delta_0$, so the equilibrium will be a node.  All trajectories entering the haline state above the nullcline will be attracted to the equilibrium along the splitting line, so there is a globally attracting equilibrium in the haline state.  Also, when $\lambda < 1$ we see that $F'(\lambda) > 0 $ if and only if $\lambda < (1+A)/(2A),$  and this case corresponds to an attracting equilibrium in the thermal state.  When the system has a global attractor, there is a fixed equilibrium circulation state.  

\begin{figure}[t]
        \centering
           \begin{subfigure}[t]{0.45\textwidth}
           	\centering
               \includegraphics[height=0.2\textheight]{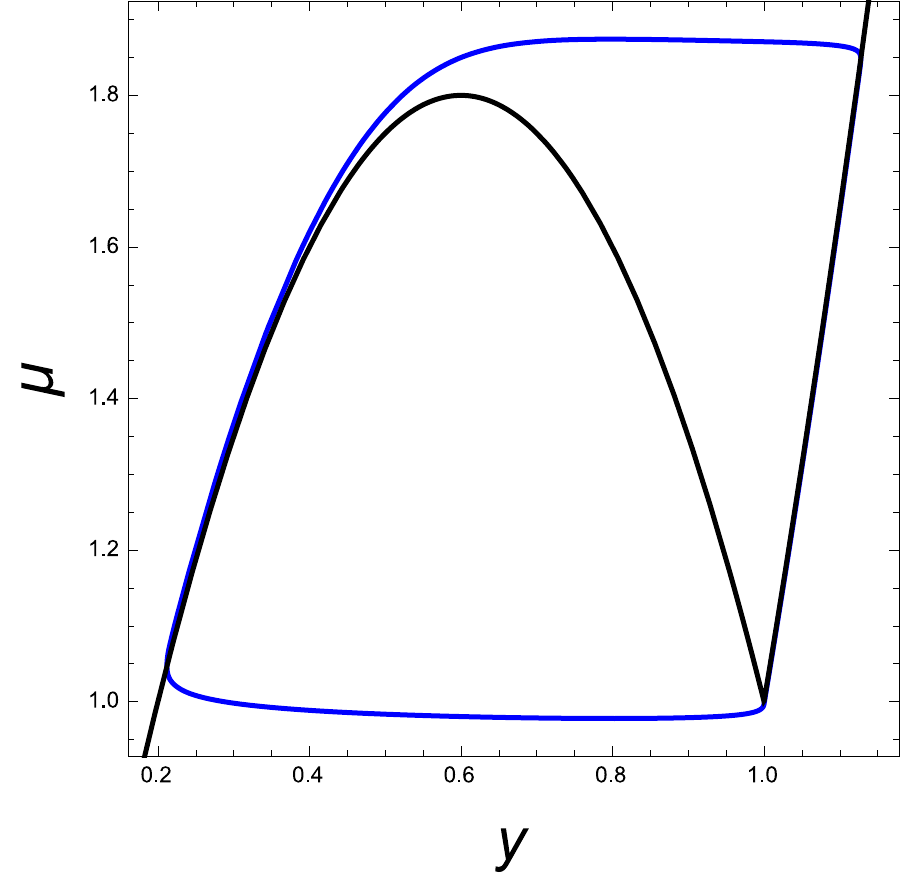}
                 \caption{Stable periodic orbit when $A=5$, $\lambda=0.8$, and $\delta=0.1$}
                \label{stomFig:ro2}
        \end{subfigure}
        ~ 
       \begin{subfigure}[t]{0.45\textwidth}
       	\centering
               \includegraphics[width=\textwidth]{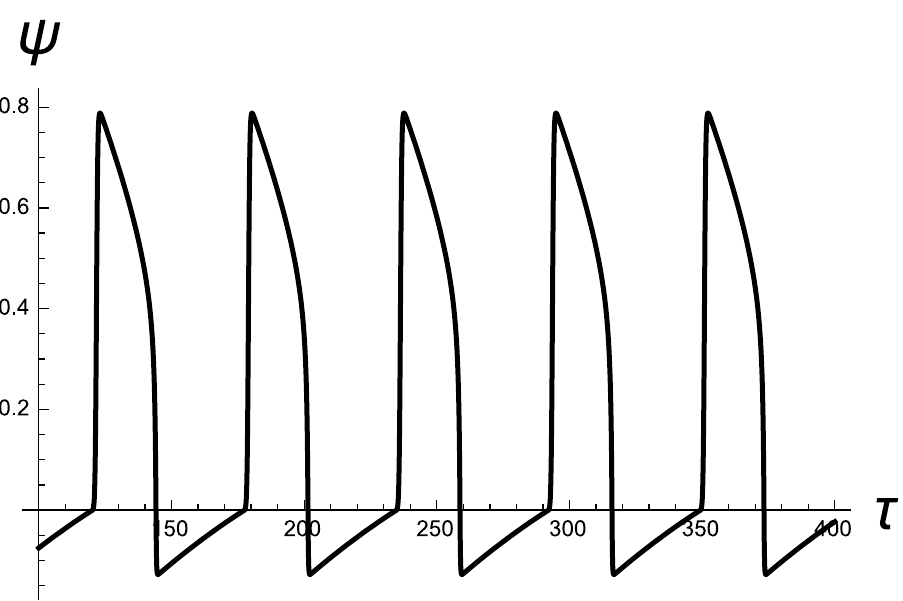}
                \caption{Time series for $\psi$ for the trajectory in (a)}
                \label{stomFig:psi2}
        \end{subfigure}
           \begin{subfigure}[t]{0.45\textwidth}
                \includegraphics[width=\textwidth]{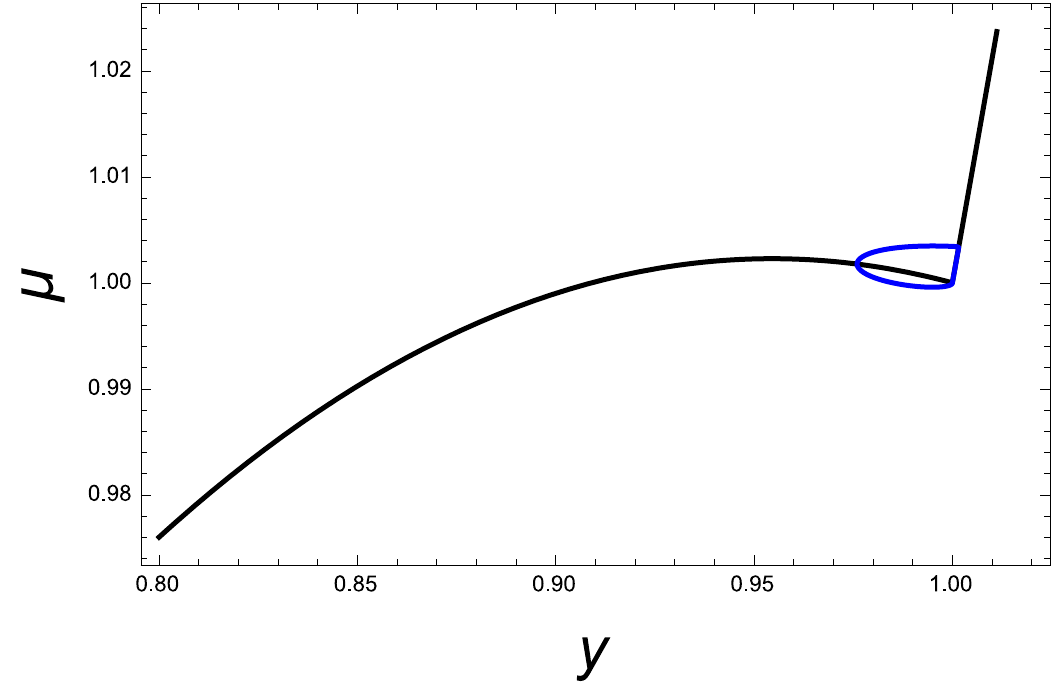}
                \caption{Canard trajectory when $A = 1.1$, $\lambda = 0.995$, and $\delta_0 = 0.01$.}
                \label{stomFig:CC}
        \end{subfigure}
        ~ 
       \begin{subfigure}[t]{0.45\textwidth}
                \includegraphics[width=\textwidth]{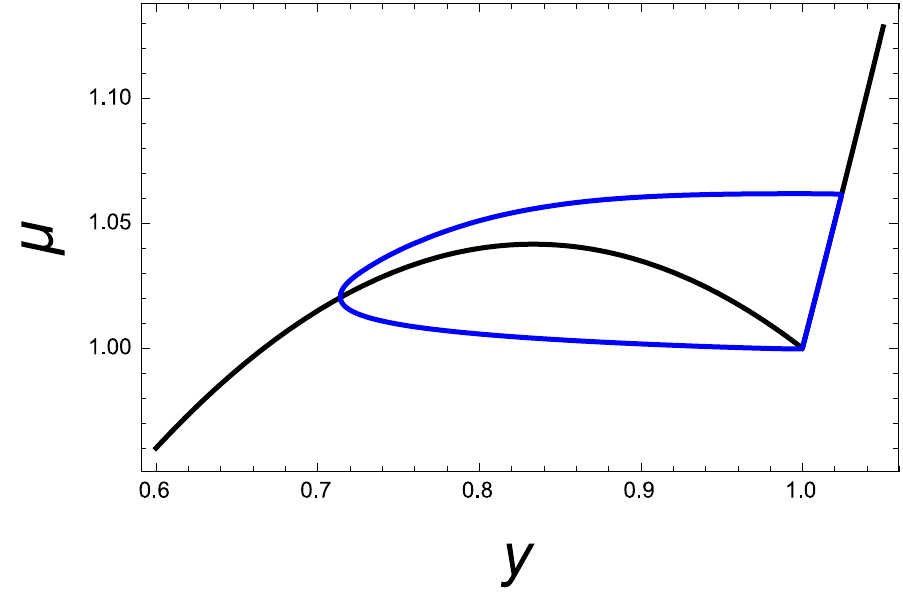}
                \caption{Super-explosion when $A=1.5$, $\lambda = 0.995$, and $\delta_0 = 0.01$.}
                \label{stomFig:SE}
        \end{subfigure}
        \caption{Oscillatory behavior in \eqref{canards}.}
        \label{stomFig:oscillations2}
\end{figure}

When $(1+A)/(2A) < \lambda < 1$, $F'(\lambda)<0$ and the equilibrium is unstable, indicating the presence of bifurcations as $\lambda$ passes through the local extrema of the critical manifold.  In particular, we discuss the nature of the bifurcation when $\lambda =1.$  If $A < 1 + 2 \sqrt{ \delta_0}$, then $F'_+(1) < 2 \sqrt{\delta_0}$ and the equilibrium will be an unstable focus near the corner.  Thus, the bifurcation will create canard cycles as shown in Figure \ref{stomFig:CC}.  However, if $A < 1 + 2 \sqrt{ \delta_0}$, then $F'_+(1) > 2 \sqrt{\delta_0}$.  The bifurcation turns a stable node into an unstable node.  In this case, the bifurcation will be a super-explosion whereby a stable relaxation orbit (bounded away from the equilibrium) appears instantaneously upon bifurcation.  The relaxation oscillation resulting from the super-explosion is depicted in Figure \ref{stomFig:SE}.  

The take-home story from this bifurcation analysis is that the phase space geometry determines the nature of the bifurcation and also the periodic orbit that the system exhibits.  Recalling Figure \ref{stomFig:o18}, the relaxation oscillations seen in the data do not have fixed amplitude.  One possible indication is that changes in the background climate state modulate the amplitude of the oscillations in circulation strength through the parameters $A$ and $\lambda$.  In the next section, we introduce orbital forcing, finding evidence that the oscillations with smallest amplitude in Figure \ref{stomFig:o18} actually correspond to canard trajectories.

\begin{remark}
The Hopf bifurcation at $\lambda = (1+A)/(2A)$ is degenerate.  When the slow nullcline intersects the critical manifold on the unstable branch, we will still have an attracting periodic orbit guaranteed by the Poincar\'{e}-Bendixson Theorem.  Simulation of \eqref{red2} with $(1+A)/(2A) <  \lambda < (1+A)/(2A) + \delta_0^2$, indicates that the system undergoes a classical canard explosion, however there also appear to be attracting periodic orbits for $\lambda < (1+A)/(2A)$.
\end{remark}		

\begin{remark}
In reality, \eqref{lin3} contains two singular perturbation parameters $\eps,\delta$, and as such represents a two-parameter relaxation oscillator.  The analysis in this section focuses on perturbing from $\delta=0$ to $\delta>0$, however we have made no mention of perturbing away from $\eps =0.$  The stability of the manifold $\{x=1\}$ indicates that our analysis will be a reasonable approximation of the dynamics for $\eps>0$ and sufficiently small.  Two-parameter singular perturbation problems have been studied in smooth systems \cite{kosiuk2011} however the full analysis of \eqref{lin3} is beyond the scope of this work.
\end{remark}

\subsection{Modulation by Orbital Forcing}

We investigate the effects of long term variations in applied forcing on the model's oscillations. During the last glacial period D-O events show significant spread in their durations, pacing and amplitude, possibly modulated by variations in the Earth's obliquity in time. The obliquity parameter controls the amplitude of the seasonal cycles and the zonal insolation gradient. Climate proxy data over several glacial cycles show significant spectral energy in the 41 kyr band \cite{wunsch2003}. According to one hypothesis the glacial terminations themselves were occurring on every second or third obliquity cycle \cite{huybers2005}. 

A variation in zonal insolation gradients will naturally affect the atmospheric temperature ($T^a$) and salinity ($S^a$) gradients. In system \eqref{canards}, the inclusion of orbital forcing implies that parameters $A$ and $\lambda$ become time-dependent.  The new system becomes
\begin{equation}
	\label{forced}
	\begin{array}{rl}
		\dot{y} =& \mu - y - A(\tau) |1-y| y \\
		\dot{\mu} =& \delta_0 ( \lambda(\tau) - y),
	\end{array}
\end{equation}
where $A(\tau)$ and $\lambda(\tau)$ are periodic functions that are scaled to $T^a(t)$ and to the ratio of the strength of temperature and salinity feedbacks respectively, and $\tau$ is the dimensionless time variable of the system.  We define the functions
\begin{align*}
	A(\tau) &= \bar{A} + p \sin \omega \tau \\
	\lambda(\tau) &= \bar{\lambda} + q \sin \omega (\tau - \theta).
\end{align*}
where $\theta$ is the phase difference between the actual atmospheric temperature gradient and the feedback response of the atmosphere due to sea ice or other climatic processes.

\begin{figure}[t]
        \centering
           \begin{subfigure}[t]{0.4\textwidth}
           	\centering
               \includegraphics[height=0.2\textheight]{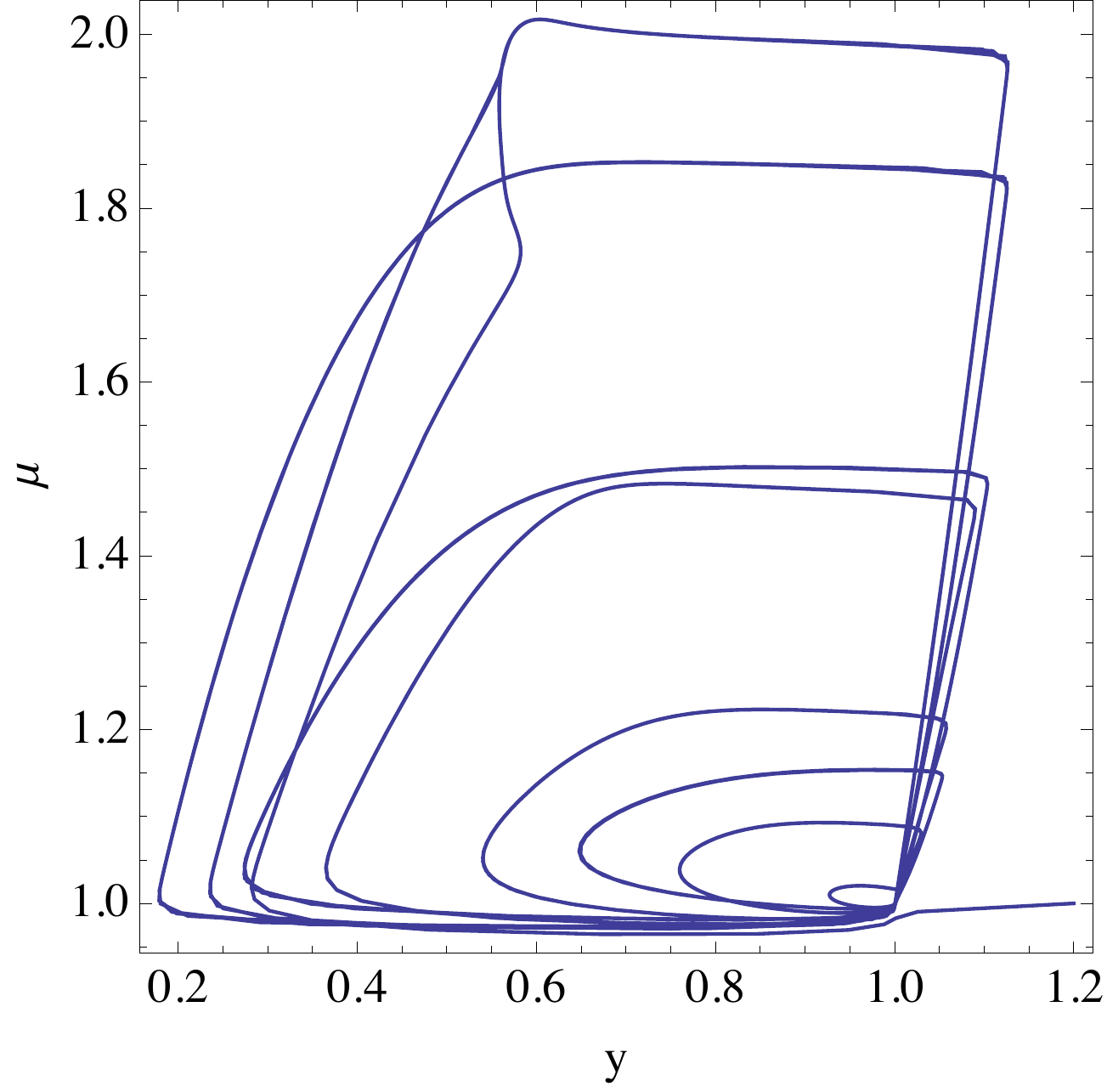}
                 \caption{Trajectory of the forced system in phase space.}
                \label{stomFig:forcedTraj}
        \end{subfigure}
        ~ 
       \begin{subfigure}[t]{0.5\textwidth}
       	\centering
               \includegraphics[width=\textwidth]{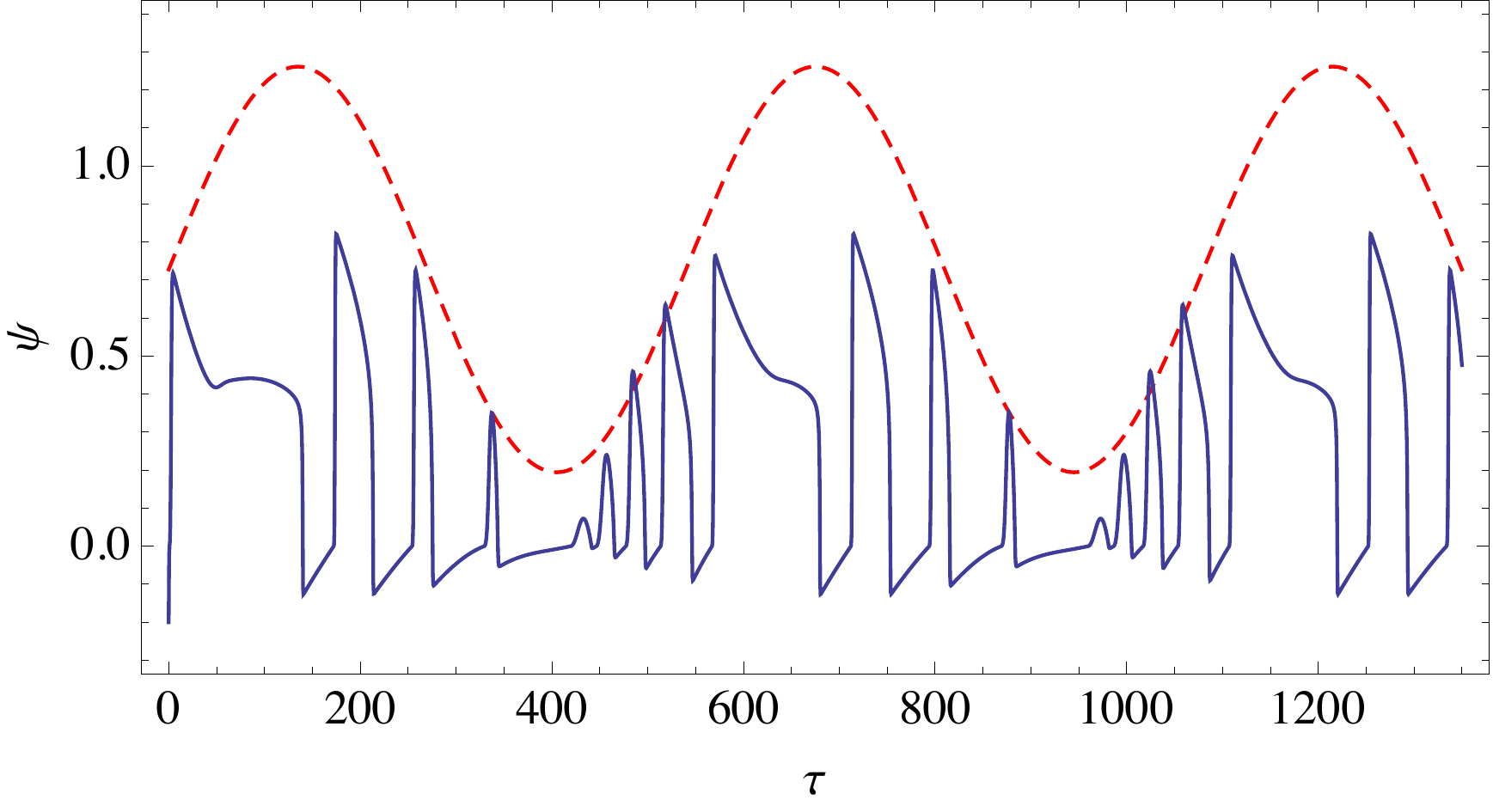}
                \caption{Time series for $\psi$ (solid) with the scaled obliquity (dashed) variations from the last 100 kyr. The units of time are arbitrary.}
                \label{stomFig:forcedPsi}
        \end{subfigure}
        \caption{Oscillations in the forced system \eqref{forced} when $\delta_0 = 0.07$, $\bar{A}=3.5$, $p = 2.4$, $\bar{\lambda} = 0.8$, $q = 1.99$, $\omega = \pi / 270,$ and $\theta = 250.$ }
        \label{stomFig:forced}
\end{figure}

We run the model with $A(\tau)$ scaled in a way so as to mimic the obliquity variation in the last 100 kyr.The model output from \eqref{forced}, shown in Figure \ref{stomFig:forcedPsi}, has important similarities to the data (Figure \ref{stomFig:o18}). Both model and data show subdued fluctuations during the dips in obliquity around 70 and 40 kyr bp, and longer interstadials during periods of maximum obliquity. The canard cycles appearing in the model output could be hidden from the temperature proxy record if the ocean atmosphere feedback is through sea ice or other processes. 

While the model output does not seek to explain every major feature of the climate record it illustrates that an intrinsic climatic oscillation could be modulated by long term variations in surface forcing. Although we do not explicitly investigate the effects of time varying salinity (freshwater) forcing, the model dynamics make it clear that it too will have  modulating effects on the underlying periodicity. The glacial climate was likely influenced by both insolation and freshwater forcing, one or the other dominating over different time windows. During the latter part of the last glacial period freshwater forcing possibly dominated due to the frequent destabilizations of the Laurentide ice sheet (Heinrich events) from cumulative ice build up on land. 

Sea surface temperature and sea ice extent proxies show that a millennial scale oscillation was also present during the warmer Holocene period that commenced at about 10 kyr before present \cite{bond1997}. These fluctuations, known as Bond events, were of diminished amplitude compared to the glacial D-O events, but had similar characteristics. They are possibly manifestations of the same oscillation mechanism responsible for the D-O events, but under a very different climate state, one where the background overturning circulation rate is higher.

\section{Discussion}

Incorporating the environmental forcing parameters as dynamic variables in Stommel's 1961 thermohaline circulation model produces a relaxation oscillator.  
In some sense, this paper is a case study for the usefulness of conceptual models.  Since the limit cycle is seen in the reduced system \eqref{red2}, the oscillation can be described by a system with two degrees of freedom, which is the minimum requirement for an oscillator.  Aside from the non-differentiability due to the absolute value term, the equations are relatively simple.  The key to generating relaxation oscillations in the model is the value of the parameter $A$, a ratio of physical parameters including the atmospheric temperature gradient.  When , $A>1$, the $y$ nullcline is `2'-shaped indicating bistability in the fast dynamics.  In reality, the nonsmooth nature of the vector field is fortuitous since the necessary bistability arises from the lack of differentiability of the absolute value function.  If the critical manifold were a smooth cubic, the system would be indistinguishable from the van der Pol system.  In fact, that is the way scientists have assumed \eqref{general3} would behave, despite GSPT breaking down at one of the jump points.  In essence, we provide rigorous justification for that assumption.

The dynamic forcing parameter $\mu$ can be thought of as feedback mechanisms between the ocean and atmospheric states. Sea ice and evaporative feedback are possible candidates that fulfill the requirements for the system to produce self-sustained oscillations. The canard cycles observed in this model can also be seen in more complicated physical models of sea ice and circulation \cite{thsaha,ppsaha}.  

Forcing the model with obliquity variations from the last 100 kyr produces important similarities with the climate record, namely colder and more subdued periods during the low points in obliquity, longer duration D-O events near peak obliquity, and the pattern of progressively declining fluctuations. In reality the climate state was possibly modulated by both orbital and freshwater forcing, where one or the other dominated during different time windows. During the latter part of the glacial period Heinrich events would have added large quantities of freshwater which would have also modulated the underlying limit cycle.

By making a philosophical shift from thinking of a hysteresis loop (i.e., a series of saddle-node bifurcations) to a relaxation oscillation (i.e., intrinsic oscillations), we are able to find a parameter regime that allows for canard trajectories.  Furthermore, it appears that these canard trajectories might play a physical role in the climate record.  While canard cycles are notoriously difficult to find in data of systems that are truly 2D, it is clear that our model \eqref{forced} is at most a backbone of the underlying model governing large-scale ocean circulation in the North Atlantic. 

\section*{Acknowledgements}
We would like to thank the Mathematics and Climate Research Network and Chris Jones for their conversations throughout the process.

\newpage

\doublespacing
\bibliographystyle{amsplain}
\bibliography{sources}

\end{document}